\documentclass[mathptm]{jcmvar}

\usepackage{graphicx}
\usepackage{amsmath}

\usepackage{xspace}

\newcounter{click}
\newcommand\ud{{\mathrm{d}}}
\newcommand\pd[2]{{\frac{\partial#1}{\partial#2}}}
\newcommand\diff[2]{{\frac{\ud#1}{\ud#2}}}
\newcommand\bP{{\mathbb{P}}}
\newcommand\bE{{\mathbb{E}}}
\newcommand\half{{\mbox{$\frac12$}}}
\newcommand\rel{\simeq}
\newcommand\shear{\gamma}

\newcommand\Sinc{{\mathsf{S}}}
\newcommand\rayt{\tau}
\newcommand\dt[1]{\diff{#1}{\rayt}}

\newcommand\xstart{x_s}
\newcommand\xend{x_e}
\newcommand\ystart{y_s}
\newcommand\yend{y_e}
\newcommand\tstart{\rayt_s}
\newcommand\tend{\rayt_e}

\newcommand\exit{diffusion-on\xspace}
\newcommand\entrance{diffusion-off\xspace}

\newlength{\figwidth}
\DeclareMathOperator{\sgn}{sgn}

\newtheorem{theorem}{Theorem}
\newtheorem{conjecture}[theorem]{Conjecture}

\begin{document}

\title[Optimally coupling the Kolmogorov diffusion]{Optimally coupling
  the Kolmogorov diffusion, and related optimal control problems}

\author{Kalvis M. Jansons}
\email{coupling@kalvis.com}
\address{Department of Mathematics,\\
University College London,\\
Gower Street,\\
London WC1E 6BT,\\ UK}

\author{Paul D. Metcalfe}
\address{Department of Applied Mathematics and Theoretical
  Physics,\\
  Centre for Mathematical Sciences,\\
  Wilberforce Road,\\ Cambridge CB3 0WA,\\ UK}

\keywords{Kolmogorov diffusion, optimal coupling, coupling
    before an exponential time, stochastic control, optimal control}

\classification{60J60, 93E20}

\begin{abstract}
  We discuss the optimal Markovian coupling before an exponential time
  of the Kolmogorov diffusion, and a class of related stochastic
  control problems in which the aim is to hit the origin before an
  exponential time.  We provide a scaling argument for the optimal
  control in the near field and use rational WKB approximation to
  obtain the optimal control in the far field, and compare these
  analytical results with numerical experiments.

  In some of these optimal control problems, in which the advection
  velocity field is bounded, we show that the probability of success
  field agrees exactly with its leading-order asymptotic approximation
  in some areas of the plane, up to an undetermined multiplicative
  constant.  We conjecture a necessary and sufficient condition for
  this behaviour, which is strongly supported by numerical
  experiments.
\end{abstract}

\maketitle

\section{Introduction}\label{sec:intro}

The \emph{coupling method} is a probabilistic tool for studying the
convergence to stationarity of a random process.  By running two
copies of the process at the same time, but imposing a dependence
between the two copies to make them collide quickly, it is possible to
obtain rigorous bounds on the convergence to stationarity of
probability distributions (in the $L_1$ norm).  This technique --- the
coupling method \cite{DurrettR:1996} --- is often used to prove
convergence to stationarity by proving that collision is guaranteed in
infinite time.

This paper discusses the optimal Markovian coupling before an
independent exponential time of two copies of the Kolmogorov
diffusion.  This study is motivated both by some recent work by
\cite{KendallWSPriceCJ:2004} and by problems in fluid-dynamical
mixing, though the latter application will be studied in depth
elsewhere.

Two copies of the Kolmogorov diffusion can be almost surely coupled in
infinite time (see \cite{BenArousGetal:1995} for proof).  To obtain a
nontrivial optimization problem, it is necessary to add some degree of
urgency to the problem.  The simplest way of doing this is by
requiring coupling before an independent exponential time; this
time-limited problem was posed by
\cite{KendallWSPriceCJ:2004}.  The memoryless nature of the
exponential time means that the resulting coupling problem remains
time-independent.  This modified coupling problem can be handled as an
optimal control problem; we use rational asymptotic methods, scaling
arguments, and numerical experiments to derive \emph{optimal}
controls, and hence optimal coupling strategies.

The fluid-dynamical motivation of the work comes from the simple model
problem of mixing in a shear flow.
The path of a passive tracer in a shear flow $v(x) = \shear x$ is
governed by the It\^o stochastic differential equations (SDE)
\begin{equation}
  \label{eq:tracer_path}
  \begin{aligned}
    \ud X_t &= \ud B_t, \\
    \ud Y_t &= \shear X_t\ud t + \ud W_t,
  \end{aligned}
\end{equation}
in which $B$ and $W$ are independent standard Brownian motions.  For
simplicity, we neglect diffusion in the streamwise direction, as is
common in dispersion theory.  Our process $(X,Y)$ then becomes the
Kolmogorov diffusion \cite{KolmogoroffA:1934}, which is the simplest physically
motivated model of thermally driven particle motion%
\footnote{i.e.\ \emph{Brown's} Brownian motion}.  This does modify the
mixing problem near the origin, but the far-field behaviour of the
mixing problem is qualitatively unchanged.  We obtain a toy
fluid-dynamical problem, but one which is of independent interest
\cite{BenArousGetal:1995,KendallWSPriceCJ:2004}.

In \S\ref{sec:formulation} we formulate our Markovian coupling problem
as an optimal control problem, and study this optimal control problem
numerically in \S\ref{sec:numerics}.  We see that the spatial
dependence of the optimal control divides the plane into geometrically
simple regions.  We identify near-field and far-field regions of the
plane and study the shapes of these geometrical features in these
asymptotic limits in \S\S\ref{sec:near_field},\ref{sec:far_field}.

Using scaling arguments, in \S\ref{sec:near_field} we study the
optimal control in the near field.  By balancing two modes of failure of
the near-field motion, we obtain a surprising scaling law for the
near-field behavior of the optimal control, which agrees well with the
numerical results.  In \S\ref{sec:far_field} we use WKB analysis
(i.e.\ the large deviation limit) to study our problem in the far
field.  We obtain a complete solution of the problem to WKB order,
which again agrees well with our numerics.  The near-field and
far-field analysis, taken together, provide a complete qualitative
picture of this optimal control problem.

We discuss some variant versions of our basic optimal control problem
in \S\ref{sec:variant}, which highlight noteworthy features of this
class of control problems that can still be understood in the WKB
limit.

Finally, in \S\ref{sec:discussion}, we give a complete description of
the optimal coupling before an exponential time of the Kolmogorov
diffusion, and discuss the extension of these ideas to fluid-dynamical
problems.

\section{Formulation}\label{sec:formulation}

The separation $(X,Y)$ of two copies of the Kolmogorov diffusion is
governed by the It\^o SDE
\begin{equation}\label{eq:separation}
  \begin{aligned}
    \ud X_t &= \ud B^{(1)}_t - \ud B^{(2)}_t, \\
    \ud Y_t &= \shear X_t \ud t,
  \end{aligned}
\end{equation}
in which $B^{(1)}$ and $B^{(2)}$ are both standard Brownian motions.
(This also governs the separation of two particles in a shear flow if
the streamwise diffusion is ignored.)
We are free to impose any
dependence between the driving Brownian motions, and we do this in such a
way as to drive the separation $(X,Y)$ to the origin quickly (see
\cite{JansonsKMMetcalfePD:2006}).

We now impose a Markovian dependence between the two driving Brownian
motions $B^{(1)}$ and $B^{(2)}$.  By a convexity argument, the local
optimization problem can be solved by requiring $\ud B^{(1)}_t = \pm \ud
B^{(2)}_t$, and our problem reduces to a control problem in two dimensions.

With little cost, we consider a richer class of control problems in
which the velocity in the $y$ direction, $v(x)$, is an
antisymmetric function of $x$ (subject to regularity requirements),
although it should be noted that only the linear case can be
interpreted as an optimal coupling problem.  We therefore choose, for
each point $(x,y)$, $\sigma(x,y)\in\{0,\sqrt{2D}\}$ in such a way as
to shepherd the process $(X,Y)$ defined by the It\^o SDE
\begin{equation}\label{eq:ito_sde_control}
  \begin{aligned}
    \ud X_t &= \sigma(X_t,Y_t)\ud B_t, \\
    \ud Y_t &= v(X_t) \ud t 
  \end{aligned}
\end{equation}
into the origin, where $B$ is a standard Brownian motion.  Our control
here is very weak; we have freedom only to switch the cross-stream
diffusion on or off, and it is not immediately obvious that we have
sufficient control to steer the particle into the origin with non-zero
probability.  It \emph{is} clear that almost all paths that eventually hit
the origin must have infinite winding number about the origin.

A simple --- degenerate --- coupling strategy allows us to bring the
path within $\epsilon$ of the origin, for arbitrary $\epsilon>0$.
Allow the particle to diffuse until it hits the line $x=\epsilon$ in
$y<0$.  The particle then hits $y=0$ at $x=\epsilon$.  Using a
sequence of such steps, the particle can then be made to hit the
origin with probability $1$ (see \cite{BenArousGetal:1995} for proof)%
\footnote{If one attempts to numerically compute a coupling strategy
  for the Kolmogorov diffusion
  that gives almost-sure coupling in infinite time, one obtains this
  degenerate strategy.  $\epsilon$ is set by the
  gridscale.}.

To remove this degeneracy it is necessary to limit the time allowed to
the particle.
A simple and natural way to do this is to mark the path
at an independent exponential time of rate $\lambda$,
and to require the particle to hit the origin before its path is marked.
This is also in keeping with our model of mixing in a fluid-dynamical
system; we require our paths to couple before they are swept apart
into the bulk of the flow.

Let $H_0$ be the (possibly infinite) time that the particle first hits
the origin, and $T_\lambda$ be the time at which the particle path is
marked.
Now, we seek to calculate
\begin{equation}
  \label{eq:1}
  \phi(x,y) \equiv \sup_{\sigma} \bP_\sigma^{(x,y)}[H_0 < T_\lambda],
\end{equation}
where $\bP^{(x,y)}_\sigma$ is the law of a particle started at
$(x,y)$, for a given control $\sigma$, and $\equiv$ denotes a definition.
Due to the antisymmetry of $v$, for all $(x,y)$,  $\phi(-x,-y) = \phi(x,y)$, and we
use this symmetry freely; in particular all our numerical solutions
are given in $x\geq 0$.

To derive the optimal control, consider a path started at
$(x,y)$.  We apply an arbitrary control $\sigma$ over the time
interval $[0,\min(h,H_0)]$, and the optimal control $\sigma^*$ over
the time interval $(\min(h,H_0),H_0]$.  Now the probability that this
path hits the origin before being marked is
\begin{equation}
  \begin{aligned}
    P &= (1 - \lambda h)\bE[\phi(X_h,Y_h)] + o(h) \\
    &= (1 - \lambda h)\left(\phi(x,y) + h v(x) \pd{\phi}{y} + \half h
      \sigma^2 \pd{^2\phi}{x^2} \right) + o(h) \\
    &= \phi(x,y) + h \left(
      v(x) \pd{\phi}{y} + \half \sigma^2 \pd{^2 \phi}{x^2} - \lambda
      \phi \right) + o(h).
  \end{aligned}
\end{equation}

We maximize this probability by choosing
\begin{equation}\label{eq:optimal_sigma}
\sigma^*(x,y) = \begin{cases}
  \sqrt{2D}, & \partial^2\phi/\partial x^2 \geq 0 \\
  0, & \mbox{otherwise},
\end{cases}
\end{equation}
and we see that $\phi$ satisfies the backward equation
\begin{equation}
  \label{eq:backward}
  v(x) \pd{\phi}{y} + D \left[\pd{^2\phi}{x^2}\right]^+ - \lambda \phi = 0,
\end{equation}
where $[X]^+ \equiv \max(0,X)$. 
Note that we require the equality in the condition of
\eqref{eq:optimal_sigma} to avoid pathologies, for example when there
is an open set in which $\partial\phi/\partial x$ is constant.

\section{Numerical results}\label{sec:numerics}

We solve the backward equation \eqref{eq:backward} by timestepping
the partial differential equation 
\begin{equation}
  \label{eq:2}
  \pd{\phi}{t} = v(x) \pd{\phi}{y} + D \left[\pd{^2\phi}{x^2}\right]^+
  - \lambda \phi,
\end{equation}
until a steady state is reached.  We impose the initial condition
\begin{equation}
  \label{eq:numerics_ic}
  \phi(0,x,y) = \begin{cases}
    1, & (x,y) = (0,0) \\
    0, & \text{otherwise},
  \end{cases}
\end{equation}
and the boundary condition $\phi(t,0,0)=1$.  Imposing the symmetry
$\phi(x,y) = \phi(-x,-y)$, we only solve on $x \ge 0$.
Using a $\tanh$ mapping,
we transform all occurrences of $(-\infty, \infty)$ to $(-1,1)$.
For unbounded problems, we use operator splitting \cite{McLachlanRIQuispelGRW:2002} to avoid CFL-based stability restrictions: on each
timestep we first
solve the hyperbolic part of \eqref{eq:2} using linear
interpolation, and second solve the diffusive part using forward
Euler in time, with second-order finite differences in $x$.  This is
unnecessary for bounded problems, where we use forward Euler
in time, first-order upwind finite differencing in $y$, and
second-order finite differencing in $x$.

\begin{figure}
  \centering
  \includegraphics[width=\figwidth]{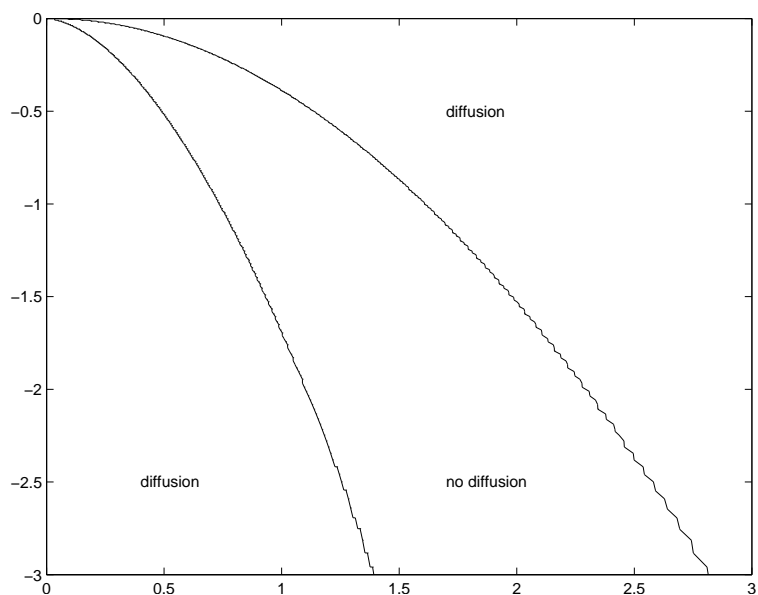}
  \caption{Boundaries of the portion of the no-diffusion region in
    $x\geq0$ for linear flow $v(x)=x$, with $\lambda = 1$, and $D=2$.
    The top boundary is $y=a(x)$ and the bottom is $y=b(x)$. The
    jagged lines for large $x$ are due to both the decreasing density
    of our stretched grid and plotting artifacts.}
  \label{fig:boundaries_linear}
\end{figure}

\begin{figure}
  \centering
  \includegraphics[width=\figwidth]{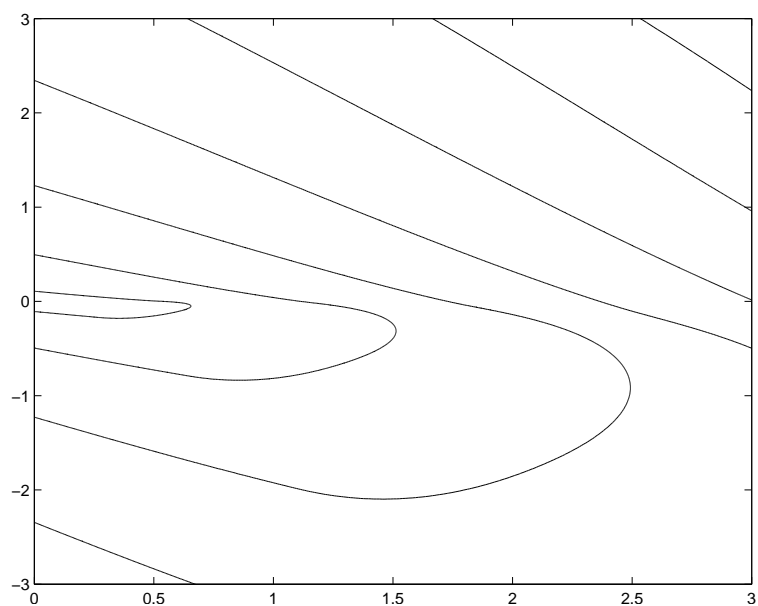}
  \caption{Contour plot of $\phi$ in $x \ge 0$ for linear flow
    $v(x)=x$, with $\lambda=1$ and $D=2$.  The contours are at
    $\phi = 3^{-n}$, for $n=1,2,\dots,7$, from around the origin outwards.}
  \label{fig:phi_field}
\end{figure}

The optimization selects a distinct region of the plane in which to
turn the diffusion off, as seen in figure \ref{fig:boundaries_linear}.
The aim of this paper is to explain the distinctive shape of this
no-diffusion region, which we sometimes refer to as region $2$.
Particles leave the no-diffusion region on the \emph{\exit boundary} $y
= a(x)$, and enter the no-diffusion region on the \emph{\entrance
boundary} $y=b(x)$.
In the far field, which we show in
\S\ref{sec:far_field} is governed by ray theory, we divide the
diffusion region into two parts: regions $1$ and $3$.  Ultimately
successful particles started in region $1$ move directly into the
near field without crossing region $2$.  Ultimately successful
particles started in region $3$ cross region $2$ in the far field
before they reach the near field.
We give a contour plot of the $\phi$ field for the case of figure
\ref{fig:boundaries_linear} in figure \ref{fig:phi_field}.

\subsection{The nature of successful paths}

A particle path that starts far from the origin, but hits the origin before
being marked, has a distinctive shape which is observed in Monte-Carlo
experiments.
In its first wind about the origin the
particle moves from the far field into the near-origin region.  This
motion is essentially ballistic, and is discussed further in
\S\ref{sec:far_field}.  After arriving in the near-origin region the
particle winds around the origin with large fluctuations,
possibly making large excursions
away from the origin.  Eventually, the particle spirals into the
origin in such a way that it never returns to the scale of a previous
spiral, and each wind about the origin results in a huge decrease in
the particle's distance from the origin.  This motion is discussed in
\S\ref{sec:near_field}.  The probability of this direct spiralling
motion increases as the particle gets closer to the origin.

Now, the probability of an unsuccessful wind, i.e.\ a wind that does
not significantly reduce the distance of the particle from the origin,
must be a rapidly decaying function of the winding number, so that
there is a non-zero probability of hitting the origin before marking
in the limit of
infinite winding number.  In particular, this implies that the
near-origin problem cannot be governed solely by Brownian scaling, as
this would give a constant (limiting) probability of failure on each
wind.  This tension against pure Brownian scaling characterizes the
dynamics of the near-origin region, and leads to the curious scaling laws
described below.

\section{Near-field scaling law}\label{sec:near_field}

Near the origin, the position of the \entrance boundary is controlled
by the balance between the two effects that lead to failure of the
directly spiralling motion: either the particle path is marked before
it hits the origin, or the particle fails to continue spiralling and
returns to the scale of its previous wind.  Because the relative
decrease in scale on each wind becomes greater and greater, the
estimate of the balance of these two sources of failure can be made
(asymptotically in small $x$, or large winding number) by considering
just one half-wind.

Now consider one half-wind for a particle near the origin starting
at a point $(x_1, y_1)$ on the \entrance boundary.  We suppose that the
flow $v(x)$ is antisymmetric, that $v(x)>0$ for $x>0$, and that $v(x)
\rel |x|^\beta$ for small $x$, where $A \rel B$ means that
$A/B$ is (at worst) a slowly varying function%
\footnote{A function $f(x)$ is slowly varying if, for all $\eta>0$, $f(\eta
  x)/f(x) \to 1$ as $x \downarrow 0$.}, which may be negative.  In this case,
the no-diffusion region lies in $xy \leq 0$.  We assume that $y_1 \rel
|x_1|^\alpha$, and we see that $x_1 y_1 < 0$.

To obtain a consistent balance, we find that the dominant contribution
to the probability of being marked on a half-wind comes from the
motion in the no-diffusion region, and has probability $F_1 \rel
y_1/v(x_1)$, for sufficiently small $x_1$.  For $v(x) \rel
|x|^\beta$, we find $F_1 \rel |x_1|^{\alpha-\beta}$.

Next, consider the particle motion from the \exit boundary until it
first hits the $y$ axis.  From Brownian scaling, the particle first
hits the $y$ axis at $y_2 \rel |x_1|^{\beta+2}$.  Now the condition
for the particle to continue to spiral into the origin is that it hits
the \entrance boundary before it makes an excursion from the $y$ axis
with $x$ displacement of scale $x_1$.  Thus this type of failure has
probability $F_2 \rel x_2/(|x_1|+|x_2|) \rel x_2/x_1$, where $x_2$ is
given by $y_2 \rel |x_2|^\alpha$, with $x_2 y_2 < 0$.  So $F_2 \rel
|x_1|^{(\beta+2)/\alpha -1}$.  The net probability of failure $F_1 +
F_2$ in this half-wind is a minimum when $F_1 \rel F_2$, which implies
that the scaling exponents are equal, giving:
\begin{equation}
\alpha^2 - \alpha (\beta-1) - (\beta+2) = 0.
\end{equation}
The scaling law thus has exponent
\begin{equation}
\label{eq:near_field_scaling}
\alpha = \half (\beta - 1 + (\beta^2 + 2\beta + 9)^{\frac12}).
\end{equation}

For the Kolmogorov diffusion, we therefore predict a scaling exponent
$\alpha = \sqrt{3}$.

We attempted to determine the slowly varying function, or possibly
constant, prefactor in the above scaling law, but were unable to
adequately control the tails of some of the distributions we
encountered (particularly in the $y$ direction).

Table \ref{tab:scaling} compares the scaling
\eqref{eq:near_field_scaling} with the results of the numerical
solution of \eqref{eq:backward}, in $-1 \leq x \leq 1$ and $-1 \leq y \leq
1$, with killing boundary conditions.  We compute the numerical
exponents by a least-squares fit of
\begin{equation}
  \label{eq:3}
  \log(-b(x)) = c + m \log x
\end{equation}
to the \entrance boundary $b(x)$, where $c$ and $m$ are constants.

\begin{table}
  \caption{Scaling exponents in the near field for various power-law flows.}
  \label{tab:scaling}
  \centering
  \begin{tabular}{c c c}
    \hline
    $\beta$ & exact $\alpha$ & numerical $\alpha$ \\
    \hline
    $1/2$ & $(\sqrt{41}-1)/4 \approx 1.351$ & $1.372$ \rule{0em}{3ex}\\
    $1$   & $\sqrt{3} \approx 1.732$ & $1.760$ \rule{0em}{3ex}\\
    $3/2$ & $(\sqrt{57}+1)/4 \approx 2.137$ & $2.201$ \rule{0em}{3ex}\\
    $2$   & $(\sqrt{17}+1)/2 \approx 2.562$ & $2.568$ \rule{0em}{3ex}
  \end{tabular}
\end{table}

This leads us to the following conjecture.
\begin{conjecture}
If $v(x)$ is antisymmetric, $v(x) > 0$ for $x > 0$, and
\begin{equation}
\lim_{x \downarrow 0} \frac{\log v(x)}{\log x} = \beta > 0,
\end{equation}
then the \entrance boundary $b(x)$ satisfies
\begin{equation}
\lim_{x \downarrow 0} \frac{\log (-b(x))}{\log x} = \half (\beta - 1 +
(\beta^2 + 2\beta + 9)^{\frac12}).
\end{equation}
\end{conjecture}

Our numerical experiments do not exclude the possibility that if $v(x)
\propto |x|^\beta$, then $b(x) \propto |x|^\alpha$, but
this clearly cannot be shown definitively by numerical means.

If, instead of the independent exponential time used elsewhere, we use
a spatially dependent marking rate $\lambda(x) \rel |x|^\Lambda$ to
limit the time, we can carry the same argument through to find
\begin{equation}
\alpha = \half(\beta-\Lambda-1 + (\beta^2 + 2\beta + 9 + \Lambda^2 - 2\beta \Lambda + 2\Lambda)^{1/2}),
\end{equation}
although it is necessary to have $\beta > \Lambda$ to obtain
a non-degenerate solution.  Numerical studies also support this
scaling exponent, but this has not been verified to as high an accuracy
as the above conjecture.

\section{WKB asymptotics in the far field}\label{sec:far_field}

As in the near field, the far-field scaling is governed by a balance
of rare events.  Here, we balance the probability that the particle path
remains unmarked with the probability that the path has a sufficient
diffusive motion to move to $x=0$.  If the particle takes a long time
to move to $x=0$, we increase the probability of the purely diffusive
motion, but we also increase the probability that the particle path is
marked, and \emph{vice versa}.

Throughout this section we assume that $v$ is antisymmetric,
continuously differentiable, and that $v'(x) > 0$ for all $x$ (where the
${}'$ denotes a derivative).  This ensures that the no-diffusion
region is connected, and its boundaries are functions of $x$.

We derive the scaling of the most important regions of the
far field by considering the $x$ and $y$ motion separately.  First, we
consider the motion in the $x$ direction.  The probability density
function for an unmarked Brownian particle in one dimension is
\begin{equation}
p(x,x_0,t) \equiv (4\pi Dt)^{-1/2} \exp\left(-\frac{(x-x_0)^2}{4Dt} - \lambda t\right),
\end{equation}
and in the limit of large $x_0>0$, the dominant contribution of unmarked
particles that hit $x=0$ is from paths of duration around $t_0$, given by
\begin{equation}
\frac{x_0^2}{4Dt_0} =  \lambda t_0,
\end{equation}
balancing the probability of achieving the large diffusive stretch
with the probability of remaining unmarked for a long time.  (This
implies a drift speed of $\sqrt{4\lambda D}$, which is exactly the
drift of a Brownian motion conditioned on tending to infinity without
being $\lambda$-marked.)  If we now consider the advection in the $y$
direction with
velocity scale $V$, 
we see that particles starting at $(x_0, y_0)$, where
\begin{equation}
y_0 \rel - V t_0
\end{equation}
will have both coordinates hit $0$ at roughly the same time.  Thus
particles starting in the neighbourhood of
\begin{equation}\label{eq:coarse_scaling}
y \rel -Vx/(\lambda D)^{1/2}
\end{equation}
make the largest contribution to the flux of unmarked particles to the
origin.  The no-diffusion region is close to this region of the
plane, so that it can control particles in this dynamically
significant region.

\subsection{The meaning of `far field'}

These ideas can be refined to give the full
leading-order behaviour of the \entrance and \exit boundaries in the
far field.  We suppose that $x$ scales as $X$, that the velocity
scales as $V$, and that $y$ scales as $X V/(4 \lambda D)^{1/2}$.
Rescaling the backward equation \eqref{eq:backward} gives
\begin{equation}
  \label{eq:4}
  \tilde{v}(\tilde{x}) \pd{\phi}{\tilde{y}}
  + \half \epsilon^2 \left[\pd{^2 \phi}{\tilde{x}^2}\right]^+ - \half
  \epsilon^{-2} \phi = 0,
\end{equation}
where tildes denote nondimensional quantities and $\epsilon^2 \equiv
(D/\lambda)^{1/2} X^{-1}$.  We now drop the tildes, and work with the
non-dimensionalized equation \eqref{eq:4} for the remainder of
\S\ref{sec:far_field}.  Next, we seek a solution of the form $\phi(x,y;\epsilon) =
\exp(-S(x,y;\epsilon)/\epsilon^2)$, and find that $S(x,y;\epsilon)$ satisfies the equation
\begin{equation}
  \label{eq:5}
  -v(x) \pd{S}{y} + \half \left[ \left(\pd{S}{x}\right)^2
    - \epsilon^2 \pd{^2S}{x^2} \right]^+ = \half,
\end{equation}
on which we impose the boundary condition $S(0,0) = 0$.  We find the
boundaries of the no-diffusion region $y=a(x)$ and $y=b(x)$ as
asymptotic expansions in $\epsilon$, as $\epsilon\to0$:
\begin{equation}
  \label{eq:6}
  \begin{aligned}
    a(x) &= a_0(x) + a_1(x) \epsilon + a_2(x) \epsilon^2 + O(\epsilon^3), \\
    b(x) &= b_0(x) + b_1(x) \epsilon + b_2(x) \epsilon^2 + O(\epsilon^3).
  \end{aligned}
\end{equation}
The expansion in $\epsilon$ is needed to split a double root.

In the diffusion regions, we seek a solution of \eqref{eq:5} in WKB
form:
\begin{equation}
  \label{eq:7}
  S(x,y;\epsilon) \sim S_0(x,y) + S_1(x,y) \epsilon + S_2(x,y) \epsilon^2 + \cdots.
\end{equation}
Terms with odd powers of $\epsilon$ cannot be ruled out in region $3$
due to $\epsilon$-scale features in region $2$.
Substituting \eqref{eq:7} into \eqref{eq:5} and expanding, we find the
equations
\begin{align}
  - v(x) \pd{S_0}{y} + \half \left( \pd{S_0}{x} \right)^2 &= \half, \label{eq:wkb_first}\\
  - v(x) \pd{S_1}{y} + \pd{S_0}{x} \pd{S_1}{x} &= 0, \label{eq:s1eq}\\
  - v(x) \pd{S_2}{y} + \pd{S_0}{x} \pd{S_2}{x} &= \half \pd{^2
    S_0}{x^2} - \half \left(\pd{S_1}{x}\right)^2 \label{eq:transport}
\end{align}
at the first three orders.

As is standard (see, for instance \cite{GoldsteinH:1980} for more details),
we solve the Hamilton--Jacobi equation \eqref{eq:wkb_first} by casting
it as its equivalent set of Hamiltonian ordinary differential equations.  We first define
a Hamiltonian
\begin{equation}
  \label{eq:8}
  H(x,y,p_x,p_y) \equiv -v(x) p_y + \half p_x^2
\end{equation}
from \eqref{eq:wkb_first}, and observe that the characteristics of
the Hamilton--Jacobi equation \eqref{eq:wkb_first} satisfy the Hamilton equations
\begin{equation}
  \label{eq:ray_eqns}
  \begin{aligned}
    \dt{x} &= \pd{H}{p_x} = p_x, \\
    \dt{y} &= \pd{H}{p_y} = - v(x), \\
    \dt{p_x} &= - \pd{H}{x} = v'(x) p_y, \\
    \dt{p_y} &= - \pd{H}{y} = 0,
  \end{aligned}
\end{equation}
where $x$, $y$, $p_x \equiv \partial{S_0}/\partial{x}$ and $p_y \equiv
\partial{S_0}/\partial{y}$ are considered as functions of the ray time
$\rayt$.  The direction of the ray time has been chosen to ensure that
information propagates away from the origin, and is in the opposite
direction to the time of particle motion. Terms such as \emph{\entrance
boundary} and \emph{\exit boundary} refer to the underlying particle
motion, and \emph{not} to the motion of rays.  We then find that the
change in $S_0$ along a ray from $(\xstart,\ystart)$ to
$(\xend,\yend)$ is
\begin{equation}
  \label{eq:9}
  \Delta S_0 \equiv \int_{(\xstart,\ystart)}^{(\xend,\yend)} \left(
    p_x \ud x + p_y \ud y \right),
\end{equation}
but that $S_1$ is constant on rays.

For a ray that remains in the diffusion regions, we introduce the quantity
\begin{equation}
  \label{eq:10}
  \Sinc(\xstart,\ystart;\xend,\yend) \equiv
\int_{(\xstart,\ystart)}^{(\xend,\yend)} \left(
    p_x \ud x + p_y \ud y \right),
\end{equation}
where the integral is taken along the ray that runs from
$(\xstart,\ystart)$ to $(\xend,\yend)$ with minimal $\Delta S_0$, so that
that  $S_0(x_1,y_1) = \Sinc(0,0;x_1,y_1)$ in region $1$.  We note that
\begin{equation}
  \label{eq:11}
  \begin{aligned}
    \pd{\Sinc}{\xend}(\xstart,\ystart;\xend,\yend) &= p_x(\tend),\\
    \pd{\Sinc}{\yend}(\xstart,\ystart;\xend,\yend) &= p_y(\tend),\\
    \pd{\Sinc}{\xstart}(\xstart,\ystart;\xend,\yend) &= -p_x(\tstart), \\
    \pd{\Sinc}{\ystart}(\xstart,\ystart;\xend,\yend) &= -p_y(\tstart),\\
  \end{aligned}
\end{equation}
where $\tstart$ and $\tend$ are the ray time at the start and end
of the ray, respectively.

We extend the idea of rays to cover the no-diffusion region; here
they have the Hamiltonian $-v(x) p_y = \half$, and we extend our definition
of $\Sinc$ with
\begin{equation}
  \label{eq:12}
  \Sinc(\xstart,\ystart;\xend,\yend) = \frac{\ystart-\yend}{2 v(\xstart)},
\end{equation}
where $(\xstart,\ystart)$ and $(\xend,\yend)$ are both in the
no-diffusion region.  Note that $\xstart=\xend$ on a ray in the
no-diffusion region.  However, $\Sinc$ notation is not used for
pieces of ray that cross from a diffusion region to a no-diffusion
region or \emph{vice versa}.

\subsection{Region 1}

In region $1$, we see that
\begin{equation}
  \label{eq:13}
  \begin{aligned}
    S_0(x,y) &= \Sinc(0,0;x,y), \\
    S_1(x,y) &= 0,
  \end{aligned}
\end{equation}
and we fit the boundary between regions $1$ and $2$ by requiring
\begin{equation}
  \label{eq:14}
  \left( \pd{S}{x} \right)^2 - \epsilon^2 \pd{^2 S}{x^2} = 0
\end{equation}
on this boundary.  This condition gives
\begin{equation}
  \label{eq:15}
  \pd{S_0}{x}(x,a_0(x)+) = 0,
\end{equation}
where the trailing $+$ denotes a limit from above.  This, not
surprisingly, is also the condition for the action to be stationary
under boundary perturbations.  Using this condition, we see that
\begin{equation}
  \label{eq:16}
  p_x = \left(1 - v(x)/v(x_0) \right)^{1/2}
\end{equation}
on a ray that hits the \exit boundary at $x=x_0$, and so
\begin{equation}
  \label{eq:17}
  a_0(x) = - \int_0^x v(\xi) 
  \left(1 - \frac{v(\xi)}{v(x)} \right)^{-1/2}\, \ud \xi.
\end{equation}

Now we can also compute the first-order correction to the \exit
boundary. On the \exit boundary, and working to $O(\epsilon^2)$, the
condition for the discriminant to change sign \eqref{eq:14} implies
\begin{equation}\label{eq:18}
\left[ \pd{^2 S_0}{x \partial y}(x,a_0(x)+) a_1(x) \right]^2 
 - \pd{^2 S_0}{x^2} (x,a_0(x)+) = 0.
\end{equation}

Next, we calculate the second-order derivatives in \eqref{eq:18}.
First, by differentiating \eqref{eq:15}, we see that
\begin{equation}
\label{eq:19}
\pd{^2 S_0}{x^2} (x,a_0(x)+) + \pd{^2 S_0}{x \partial y} (x,a_0(x)+) a_0'(x) = 0.
\end{equation}

Since the rays are vertical just above the \exit boundary ($p_x=0$,
and $p_y = -(2 v(x))^{-1}$), we
compute $\partial^2 S_0/\partial x\partial y$ just above the \exit boundary by
differentiating along rays:
\begin{equation}\label{eq:20}
\pd{^2 S_0}{x \partial y}(x,a_0(x)+) 
= \frac{\ud p_x}{\ud\rayt} \left( \frac{\ud y}{\ud\rayt} \right)^{-1}
= - v'(x) p_y/v(x) = - \left(\frac{1}{2 v(x)}\right)',
\end{equation}
and so
\begin{equation}
  \label{eq:21}
  \pd{^2 S_0}{x^2} (x,a_0(x)+) = a_0'(x) \left(\frac{1}{2
    v(x)}\right)'.
\end{equation}

Now, we use the second derivatives (\ref{eq:20}, \ref{eq:21}) in the
perturbation discriminant \eqref{eq:18} to find the
first-order correction to the \exit boundary:
\begin{equation}
  \label{eq:22}
  a_1^2(x)
  = \pd{^2 S_0}{x^2} (x,a_0(x)+) 
  \left( \pd{^2 S_0}{x \partial y}(x,a_0(x)+) \right)^{-2}
  = a_0'(x) \left[ \left(\frac{1}{2 v(x)}\right)' \right]^{-1}.
\end{equation}
We choose the positive root of \eqref{eq:22} to keep
the discriminant positive in region $1$.

Projecting our asymptotic expansion for $a(x)$ into the class of
decreasing functions, we obtain the approximation
\begin{equation}
  \label{eq:23}
  \bar{a}(x) \equiv a_0(x) \left(1-\epsilon \frac{a_1(x)}{a_0(x)} \right)^{-1},
\end{equation}
so that $a(x)-\bar{a}(x) = O(\epsilon^2)$.

The action $S_0(x,a_0(x))$ is then
\begin{equation}
  \label{eq:24}
  S_0(x,a_0(x)) = - \int_0^x \frac{a_0'(\xi)}{2 v(\xi)}\, \ud \xi.
\end{equation}

We find, then, that for the power-law flow $v(x) = \sgn(x) |x|^\beta$,
the first two terms in the asymptotic expansion of the
\entrance boundary $a(x)$ in $x>0$ are
\begin{equation}
  \label{eq:25}
  \begin{aligned}
    a_0(x) &= - x^{\beta+1} B(\half, 1+\beta^{-1})/\beta \\
    a_1(x) &= x^{\beta+1/2} \left[2 (\beta^{-1}+\beta^{-2})
      B(\half, 1+\beta^{-1})\right]^{1/2},
  \end{aligned}
\end{equation}
and, just for interest, that
\begin{equation}
  \label{eq:26}
  S_0(x,a_0(x)) = x (1+\beta^{-1}) B(\half,1+\beta^{-1})/2,
\end{equation}
where $B$ is the Beta function.

\subsection{Region 2}

To find the boundary $b(x)$ of region $2$ we find the points at which
the discriminant \eqref{eq:14} is zero.  The action in region $2$ is
easily computed:
\begin{equation}
  \label{eq:27}
  \begin{aligned}
    S(x,y) &= S(x,a(x)) + \frac{a(x)-y}{2 v(x)} \\
    &= S_0(x,a_0(x)) + \frac{a_0(x) - y}{2 v(x)} + O(\epsilon^2).
  \end{aligned}
\end{equation}
Substituting \eqref{eq:27} into the discriminant, we find $b_0 = a_0$,
and then at first order we find
\begin{equation}
  \label{eq:28}
  b_1^2(x) = a_0'(x) \left[ \left(\frac{1}{2 v(x)}\right)' \right]^{-1}.
\end{equation}
Here, we must choose the negative root to keep the discriminant
negative in region $2$; we see that $b_1 = -a_1$.  We also see that
\eqref{eq:28} has a root that gives the perturbation to the \exit
boundary; this is reassuring, as it shows that the analysis in region
$1$ is consistent with the analysis in region $2$.  Note that the
symmetry between the \entrance and \exit boundaries about $y =
a_0(x)$ is not expected to extend to $a_2$ and $b_2$.

\subsection{Region 3}

The main purpose of \S\ref{sec:far_field} is to determine the
far-field approximations to the \entrance and \exit boundaries, which
can be done without reference to region $3$.  However, there are still
some interesting features of the field in region $3$, which are
discussed here.

By direct calculation, with error $O(\epsilon^2)$, the action in
region $2$ is the analytic continuation of the action in region $1$.
Hence, the action in region $3$ is, with error $O(\epsilon^2)$, the
analytic continuation of the action in region $1$.  Now the ray tubes
in region 2 are perturbed (in position and cross-section) by
$O(\epsilon^2)$ due to the lack of diffusion, and they spend only
$O(\epsilon)$ ray time, $\rayt$, in region $2$.  This suggests that the
true difference between the analytic continuation of the region $1$
solution to region $3$ and the true region $3$ solution could be
as small as $O(\epsilon^3)$.

\begin{figure}
  \centering
  \includegraphics[width=\figwidth]{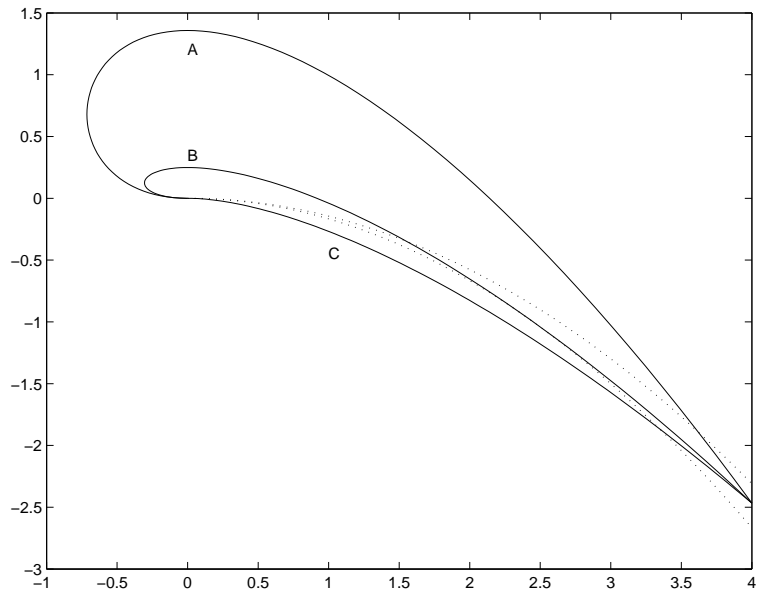}
  \caption{Rays from the origin to a point just inside region $1$, in
    linear flow (but ignoring the small effect from the no-diffusion
    region in $x<0$, $y>0$).  The rays are plotted as solid lines, and
    the boundary between regions $1$ and $3$, and the line in region
    $1$ on which non-optimal rays appear, are plotted as dotted lines.
    The caustic appears when the action along ray \textsf{A} and ray
    \textsf{C} are equal.  Ray \textsf{B} is a local \emph{maximum} of
    the action.}
  \label{fig:raypic}
\end{figure}

The only remaining feature of region $3$ is its boundary $y=c(x)$ with region
$1$.  This boundary is a caustic of the Hamilton--Jacobi problem
\eqref{eq:wkb_first}.  In region $1$, optimum rays head straight away
from the origin, moving straight into positive $x$, whereas in region
$3$ optimum rays move into negative $x$ before looping back into
positive $x$.  There is a transition region in the neighbourhood of
the boundary of regions $1$ and $3$ where non-optimal rays exist;
these rays are plotted for linear flow ($v(x)=x$) in figure
\ref{fig:raypic}, in which we launch rays from the origin to a point
just inside the boundary of region $1$.  There are three rays from the
origin to this point.  Ray \textsf{A} first moves into $x<0$, and has
a large loop there before reaching $x=0$.  Ray \textsf{B} also moves
into $x<0$, but only has a small tight loop there.  Finally, ray
\textsf{C} moves directly into $x>0$.  Rays \textsf{A} and \textsf{C}
give local minima of the action, but ray \textsf{B} gives a local
maximum.  By definition, the boundary $y=c(x)$ between regions $1$ and $3$ is
the curve on which the actions of rays of type \textsf{A} and type
\textsf{C} are equal.  Although figure \ref{fig:raypic} only shows the
behaviour of rays in linear flow, we believe the qualitative behaviour
of these rays is generic.

For linear flow we explicitly determine the action in $x\geq 0$, and
can therefore explicitly find this caustic.  Define
\begin{equation}
  \label{eq:29}
  S_\pm(x,y) \equiv \mbox{$\frac{2}{3}$} (x^2/\tau_\pm(x,y) + \tau_\pm(x,y)) 
  \pm  \mbox{$\frac{1}{3}$} x
\end{equation}
where
\begin{equation}
  \label{eq:30}
  \tau_\pm(x,y) \equiv (x^2 \pm 6 y)^{1/2} \pm x,
\end{equation}
and then $S(x,y) = S_-(x,y) + O(\epsilon^2)$ in $y\leq c(x)$, and
$S(x,y) = S_+(x,y) + O(\epsilon^2)$ in $y\geq c(x)$.  We find the
caustic $y=c(x)$ by solving $S_+(x,y) = S_-(x,y) + O(\epsilon^2)$,
giving
\begin{equation}
  \label{eq:31}
  c(x) =  - x^2/(4 \sqrt{3}) + O(\epsilon^2).
\end{equation}

\subsection{Far-field discussion}

\begin{figure}
  \centering
  \includegraphics[width=\figwidth]{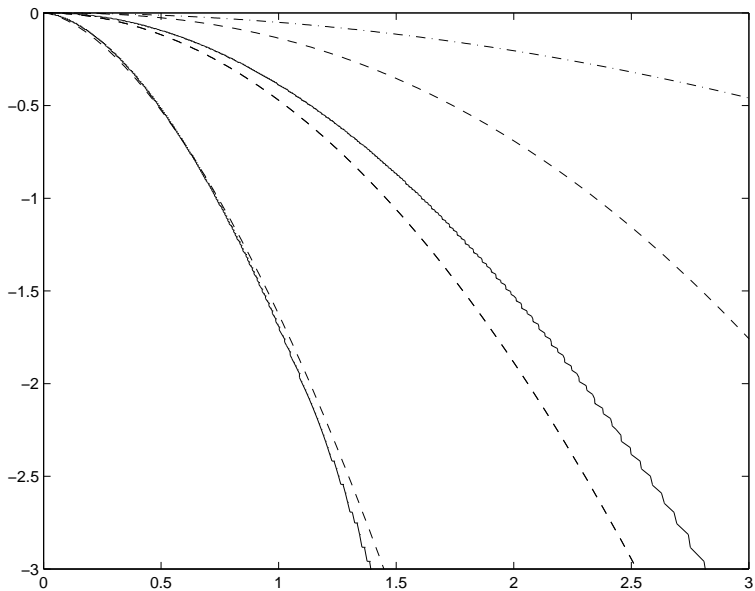}
  \caption{Comparison of numerical results and far-field asymptotics
    for the boundaries of the no-diffusion region, for linear flow
    $v(x)=x$, with $\lambda = 1$, and $D=2$.  From top to bottom, the
    numerical \exit and \entrance boundaries are plotted with solid
    lines, and the dimensional forms of the curves $y=\bar{a}(x)$,
    $y=a_0(x)$ and our asymptotic approximation for $y=b(x)$ are
    plotted with dashed lines.  The boundary, $y=c(x)$, between
    regions $1$ and $3$ is plotted with a dash-dotted line.  The
    agreement between asymptotics and numerics is reasonable, given
    that $\epsilon \approx 0.8$, based on $X=2$.}
  \label{fig:wkb_comparison}
\end{figure}

Figure \ref{fig:wkb_comparison} plots the asymptotic results of
\S\ref{sec:far_field} for the far-field \entrance and \exit boundaries,
and gives a comparison with the numerically computed boundaries, for
the case of a linear flow.  Despite the large asymptotic parameter
($\epsilon \approx 0.8$), there is very good agreement between the
asymptotic approximation of the \entrance boundary and the numerical
boundary.  The fit of the \exit boundary is worse, although there
is still reasonable agreement between the leading-order approximation
$a_0(x)$ and the numerical boundary, which is between $a_0(x)$ and its
first correction $\bar{a}(x)$.

This large deviation calculation may equivalently be done in
terms of a path integral formulation of the problem, in which the
ray equations \eqref{eq:ray_eqns} arise as the equations
governing the \emph{classical path}.

\section{Variant problems}\label{sec:variant}

\subsection{Exact snap to far-field asymptotics}

\begin{figure}
  \centering
  \includegraphics[width=\figwidth]{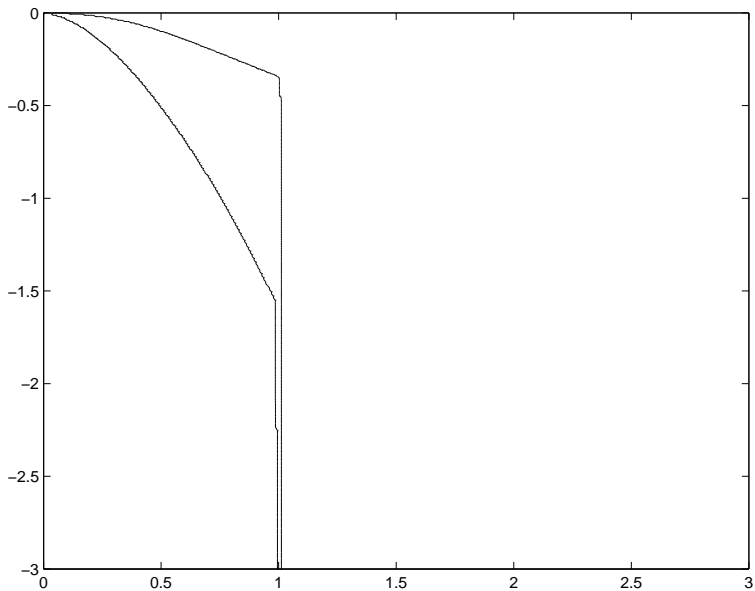}
  \caption{A no-diffusion region degenerating to a line, for $D=2$,
    $\lambda=1$, and with the flow field \eqref{eq:exact_snap_flow}.}
  \label{fig:exact_snap}
\end{figure}

Throughout this subsection, we assume that $v$ is antisymmetric,
continuous, $v(x) > 0$ for $x>0$, and without loss of generality that
$\sup_{x}v(x) = 1$.  We further assume that there is a finite $x_m
\equiv \inf \{x : v(x) = 1 \}$. Under a local `sharpness' condition of $v$ in
the neighbourhood of $x_m$ (given below), for sufficiently large $|y|$
the no-diffusion region degenerates to the lines $x= - \sgn(y) x_m$, on
which the field decays exponentially in $|y|$.  This is seen in figure
\ref{fig:exact_snap}, which shows this phenomenon with the flow field
\begin{equation}
  \label{eq:exact_snap_flow}
  v(x) = \begin{cases}
    x, & |x| < 1 \\
    \sgn x, & \mbox{otherwise}.
  \end{cases}
\end{equation}

Based on our WKB analysis, we conjecture a necessary and sufficient
condition for the no-diffusion region to degenerate to lines
$x=-\sgn(y)x_m$ in sufficiently large $|y|$.  We simply ask whether a
ray from the origin hits $x=x_m$ with $p_x=0$ at a finite ray time
(which is the condition for the ray to enter the no-diffusion
region). This reduces to the condition
\begin{equation}
\label{eq:wkb_condition}
\int_0^{x_m} \left(1 - \frac{v(x)}{v(x_m)}\right)^{-1/2} \ud x < \infty,
\end{equation}
for the no-diffusion region to degenerate to a line, and not
otherwise.  Based on numerical tests the above condition also appears
valid for the case of bounded $x$, with reflecting boundary
conditions, although the WKB argument given above is not directly
applicable because of caustics due to boundary reflections.

An example for which the above condition predicts no snap is
$v(x) = \max(x(2-x),1)$ for $x\geq0$, and this was numerically
verified. This example is close to the marginal case.

In some systems, and in some regions of the plane, this snap of the
no-diffusion region leads to a curious behaviour in which the field
$\phi(x,y)$ agrees exactly with its leading-order far-field asymptotic
solution, up to an undetermined multiplicative constant.
We refer to this phenomenon as an \emph{exact snap}.

We now consider this exact snap in detail.  From \eqref{eq:4},
$\phi(x,y) \leq \exp(-\half\epsilon^{-2} |y|)$ for all $(x,y)$, based
on the maximum advection velocity.  Define $y_m$ as the infimum of
$|y|$ in the degenerate part of the no-diffusion region.

Define $\psi$ so that
\begin{equation}
\label{eq:snap1}
\phi(x,y) \equiv \exp(-\half\epsilon^{-2} |y|) \psi(x,y),
\end{equation}
then equation \ref{eq:4} implies for $|y|>y_m$ that
\begin{equation}
\label{eq:snap2}
v(x) \pd\psi{y} + \half \epsilon^2 \pd{^2 \psi}{x^2} 
- \half \epsilon^{-2} (1 + \sgn(y) v(x)) \psi = 0,
\end{equation}
with discontinuities in $\partial\psi/\partial x$ at $x = \pm x_m$ in
respectively $y<-y_m$ and $y>y_m$.

Under the above condition for a degenerate no-diffusion region, we
find areas of the plane for which $\partial\psi/\partial y = 0$,
namely $y < -y_m$ and $x < x_m$, and the corresponding area in $y>0$.
Equation \eqref{eq:snap2} then simplifies, and we are able to find the
field in these areas of the plane, up to an undetermined constant
$\phi(x_m,y_m)$.

For $v(x) = \max(x,1)$, for $x\geq 0$, we see that $x_m=1$, and that our
WKB condition \eqref{eq:wkb_condition} is satisfied.  In our exact
snap areas, \eqref{eq:snap2} can be exactly solved in terms of Airy
functions and exponentials.

We test the hypothesis that $\partial \psi/\partial y = 0$ in $x \leq 1$
for sufficiently large $|y|$ by computing
\begin{equation}
  \label{eq:32}
  {\cal N}(\xi_1, \xi_2; y_1,y_2) 
  \equiv \sup_{\xi_1\leq x\leq\xi_2} \frac{|\psi(x,y_1) -
    \psi(x,y_2)|}{
    \psi(x,y_1) + \psi(x,y_2)}.
\end{equation}

For the data of figure \ref{fig:exact_snap}, we find that ${\cal
  N}(-\infty, 1;-5.50,-7.46) = 5 \times 10^{-4}$.
(In figure \ref{fig:exact_snap}, these dimensionless values
  of $y$ correspond to dimensional values of $-1.95$ and $-2.64$
  respectively, based on $X=1$ and $V=1$.  This gives
  $\epsilon^4=2$.)
It is difficult to get stronger numerical evidence for this in a
doubly infinite system; our stretched grid coarsens away from the
origin, and we choose $y$ stations far enough from the ends of the
grid to ensure that the field is reasonably well resolved.  For the
values of $y$ chosen here, the grid spacing in $y$ is approximately
$0.06$.  We find, from our numerical results, that the decay of the
field along the filament of the no-diffusion region on
$x=-\sgn(y) x_m$, $|y| > y_m$ is $\exp (-0.998 \times \half\epsilon^{-2}
y)$.
Note, however, that
the correlation coefficient between the grid points $(x_i)$ and
$(\log\phi(x_i,-5.50))$, in $x\leq -1$, is $1-10^{-8}$.
(We work in $x \leq -1$ because, for this problem, the predicted
  $\psi$ field is particularly simple here.) 
The exponential $(C_1 \exp (\sqrt{2} \epsilon^{-2} x_i))$ (where $C_1$ is
a constant) is an excellent fit to the $\phi$ field on $x<-1$,
$y=-5.50$.  We also find that $C_2 \exp (\sqrt{2} \epsilon^{-2}x +
\half\epsilon^{-2} y)$ (where $C_2$ is a constant) is a good fit
to the $\phi$ field over the area $y<-y_m$, $x<-1$.  However, the
decay at fixed $y<-y_m$ of $\phi(x,y)$ in $x>1$ is
(numerically) far from $\exp(-\sqrt{2} \epsilon^{-2} x)$, demonstrating
that there is no exact snap in this area.

A more convincing demonstration of this exact snap occurs in
systems with bounded $x$ and reflecting boundary conditions.  In the
corresponding system with $x \in [-1,1]$ and unbounded $y$, higher
resolution in $y$ is possible, giving for example ${\cal N}(-1,
1;-5.49,-6.62)$ $< 10^{-11}$.  Here, the no-diffusion region degenerates
onto $x= - \sgn(y)$ for sufficiently large $|y|$, and we are able to
accurately resolve this boundary.

The reason we believe that some systems have an exact snap is that
this is the slowest decaying mode, so if the optimization can match to
this solution it can do no better than this for larger $|y|$, given
the amplitude of $\phi$ at the matching $y$ value.  It appears that an
exact snap does not occur in any area of the plane for which there
exists a possible particle path to the $x$ axis that does not pass
through the no-diffusion region.

\subsection{Disconnected no-diffusion regions}

In the WKB calculation of \S\ref{sec:far_field} we assume that $v$ is
antisymmetric, continuously differentiable with $v'(x)> 0$, 
the last condition ensuring that the no-diffusion region is connected,
significantly simplifying the calculation.  If however we allow
intervals $|x| \in (\xi, \eta)$ where $v'(x) = 0$, these lead to gaps in
the no-diffusion region over a range which includes $(\xi,\eta)$.

\begin{figure}
  \centering
  \includegraphics[width=\figwidth]{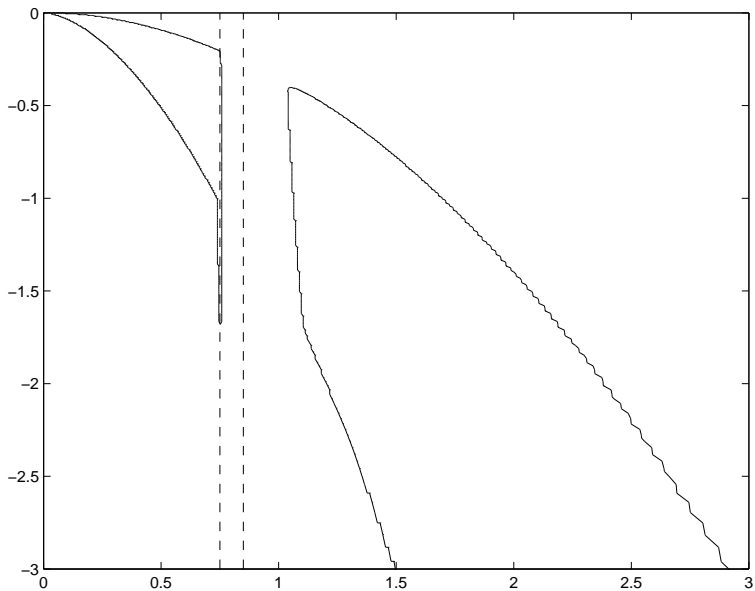}
  \caption{A gap in the no-diffusion region, for the flow field
    \eqref{eq:gap_flow}, with $D=2$ and $\lambda=1$.  The dashed lines
    show the ends of the interval $(0.75,0.85)$ in which $v$ is
    constant.}
  \label{fig:gap}
\end{figure}

An example for $v$ with
\begin{equation}
  \label{eq:gap_flow}
  v'(x) = \begin{cases}
    0, & |x| \in (0.75, 0.85)\\
    1, & \mbox{otherwise}
  \end{cases}
\end{equation}
is given in figure \ref{fig:gap}, clearly showing a gap in the
no-diffusion region due to the interval in which $v$ is flat.  Notice
that the width of the gap in the no-diffusion region is larger than
the width of the interval in which the velocity is flat, and this is a
generic feature of such systems.  Also generic is the filament of
no-diffusion region that hangs down from the part of the no-diffusion
region connected to the origin.  This filament shortens as the gap
size is reduced, and the $y$ value at which the filament ends always
roughly matches the $y$ value at the bottom of the other side of the
gap.

If there are many such small gaps close together, the no-diffusion
region is vertically striped, but its behaviour is almost the same as
that of its envelope.  This is because when a particle moves into the
envelope, it either (i) enters a very $x$-restricted diffusion region
which it passes through, (ii) enters a very $x$-restricted diffusion
region and later hits the adjacent no-diffusion region, or (iii) goes
straight into the no-diffusion region.  In all cases, the particle
leaves the envelope with almost the same $x$.

An extreme example of a striped envelope is when, apart from in the
neighbourhood of the origin needed for success, the velocity is piecewise
constant. In this case the no-diffusion strips are of zero width.

\section{Discussion}\label{sec:discussion}

We are now in a position to clearly describe the behaviour of a
particle that ultimately hits the origin before its path is
$\lambda$-marked.
For definiteness, we consider a particle in the unbounded
shear flow $v(x)=x$, started at $(0,y_0)$ with $y_0<0$ and well into
the far field.  The particle first moves out into large $x>0$, to
increase its flow speed towards the $x$ axis.  When the particle
reaches $y=y_0/2$ it turns around, and heads back towards $x=0$.  Near
this maximum of $x$, the particle gains a small advantage by turning
diffusion off, and hence truncating a quadratic
maximum in the $x$ motion, as seen in \S\ref{sec:far_field}.  The far-field
motion (in the WKB limit) is little affected by the no-diffusion
region, with corrections to the rays and the action first appearing at
$O(\epsilon^2)$ (and possibly higher), whereas the no-diffusion
boundaries have corrections that scale as $\epsilon$, caused by the splitting
of a double root.  The structure of the no-diffusion region for
$\epsilon=1$ is very similar to its structure in the far-field limit
(i.e.\ small $\epsilon$).

When the particle finally arrives near $x=0$, it winds around the
origin, with large fluctuations, until eventually it spirals directly into the origin.
Controlling this rapid spiralling motion fixes the behaviour of the
\entrance boundary near the origin, as seen in \S\ref{sec:near_field}.
In the near field the no-diffusion region is essential, as without it
the probability of hitting the origin is zero.  The near-field scaling
laws are supported well by numerical experiments (including cases with
power-law killing).  A natural question is whether such scaling laws
are possible for higher-order iterated diffusions, and whether it is
possible to go beyond scaling laws to asymptotic expansions in the
near field.

The no-diffusion region is intrinsically governed by a balance between
the probability that the particle path is marked, and the probability that a
particle hits the origin at all.  Considering the optimum no-diffusion
region as a function of $\lambda$, we see that as $\lambda$ increases,
the probability of a particle starting at a fixed point hitting the
origin marked or unmarked decreases.

The natural fluid-dynamical extension of this work is to optimal
coupling in a linear flow.  The separation $X$ between the two particles
is governed by the It\^o SDE
\begin{equation}
  \label{eq:svd_sde}
  \ud X_t = A \cdot X_t\, \ud t + (I - R(X_t)) \cdot \ud B_t,
\end{equation}
where $B$ is Brownian motion, $A$ is a traceless constant matrix, $I$
is the identity matrix, and $R$ is an orthogonal matrix.
For each point, $x$, we choose $R(x)$ in order to
rapidly shepherd $X$ into the origin, and hence achieve rapid coupling.
We do this by forming the backward equation for the probability of
successful coupling, analogous to \eqref{eq:backward}.  The local
optimization problem, in this paper only a discrete choice, now
becomes a continuum optimization problem, which may be solved by
singular value decomposition \cite{JansonsKMMetcalfePD:2006}.
The far field is still governed by ray theory, and similar techniques
to those of \S\ref{sec:far_field} will still apply.  However, the
region corresponding to the no-diffusion region of this paper will
have more structure, because the optimization is no longer a discrete
choice.
It may be possible to produce systems of the form \eqref{eq:svd_sde}
in which there is an intermediate scale resembling the near-field
limit of \S\ref{sec:near_field}.  However, it is not clear if this
effect will be important in real fluid-dynamical applications.

\bibliography{coupling}

\end{document}